\documentclass[reqno]{amsart}
\usepackage{amsthm,amsmath,amsfonts,amssymb,enumerate,latexsym,mathrsfs}

\usepackage[dvips]{graphicx}
\newtheorem{thm}{Theorem}

\newtheorem*{thm*}{Theorem}

\newtheorem*{defn*}{Definition}
\newtheorem{prop}{Proposition}

\newtheorem*{claim}{Claim}

\newtheorem{lemma}{Lemma}

\newtheorem*{cor*}{Corollary}

\newcommand{\R}{\mathbb{R}}
\newcommand{\Q}{\mathbb{Q}}

\newcommand{\funct}[2]{#1 \longrightarrow #2}

\newcommand{\Ur}{\textbf{U}}

\newcommand{\s}{\textbf{S}}

\newcommand{\dom}{\mbox{$\mathrm{dom}$}}
\newcommand{\ran}{\mbox{$\mathrm{ran}$}}
\newcommand{\vep}{\varepsilon}

\newcommand{\m}[1]{\textbf{#1}}

\newcommand{\mc}[1]{\widetilde{\textbf{#1}}}

\author{J. Lopez-Abad and L. Nguyen Van Th\'{e}}

\address{Equipe de Logique Math\'{e}matique, UFR de Math\'{e}matiques,
 (case 7012), Universit\'{e} Denis Diderot Paris 7, 2 Place Jussieu,
 75251 Paris Cedex 05, France.}

\email{abad@logique.jussieu.fr}

\email{nguyenl@logique.jussieu.fr}

\title{The oscillation stability problem for the Urysohn sphere: A combinatorial approach.}


\date{June, 2007}

\begin{document}

\begin{abstract}

We study the oscillation stability problem for the Urysohn sphere, an analog of  the distortion
problem for $\ell _2$ in the context of the Urysohn space $\Ur$. In particular, we show that this
problem reduces to a purely combinatorial problem involving a family of countable ultrahomogeneous
metric spaces with finitely many distances.

\end{abstract}

\maketitle

\section{Introduction.}

\label{section:Definitions and notations}

The purpose of this note is to present several partial results  related to the \emph{oscillation
stability problem for the Urysohn sphere}, a problem about to the geometry of the Urysohn space
$\Ur$ which can, in some sense, be seen as an analog for $\Ur$ of the well-known \emph{distortion
problem for $\ell _2$}. This latter problem appeared after the following central result in
geometric functional analysis established by Milman: For $N \in \omega$ strictly positive, let
$\mathbb{S}^N$ denote the unit sphere of the $(N+1)$-dimensional Euclidean space and let
$\mathbb{S} ^{\infty}$ denote the unit sphere of the Hilbert space $\ell _2$. If $\m{X} = (X, d^{\m{X}})$ is
a metric space, $Y \subset X$ and $\varepsilon > 0$, let also
\[(Y) _{\varepsilon} = \{ x \in X : \exists y \in Y \ \ d ^{\m{X}} (x,y) \leqslant \varepsilon \}.\]
Then:

\begin{thm*}[Milman \cite{Mil}]
\label{thm:Milman'} Let $\gamma$ be a finite partition of $\mathbb{S} ^{\infty}$. Then for every
$\varepsilon > 0$ and every  $N \in \omega$, there is $A \in \gamma$ and an isometric copy
$\widetilde{\mathbb{S}}^N$ of $\mathbb{S}^N$ in $\mathbb{S} ^{\infty}$ such that
$\widetilde{\mathbb{S}}^N \subset (A)_{\varepsilon}$.
\end{thm*}
Whether Milman's theorem still holds when $N$ is replaced by $\infty$ is the content of the
distortion problem for $\ell _2$. Equivalently, if $\varepsilon > 0$ and $f : \funct{\mathbb{S}^{\infty}}{\R}$ is bounded and
uniformly continuous, is there a closed infinite-dimensional subspace $V$ of $\ell
_2$ such that $\sup \{ \left| f(x) - f(y) \right| : x, y \in V \cap \mathbb{S} ^{\infty} \} < \varepsilon$? This
question remained unanswered for about 30 years, until the solution of Odell and Schlumprecht in
\cite{OS}:

\begin{thm*}[Odell-Schlumprecht \cite{OS}]

\label{thm:Odell-Schlumprecht}

There is a finite partition $\gamma$ of $\mathbb{S} ^{\infty}$ and  $\varepsilon > 0$ such that no
$(A)_{\varepsilon}$ for $A \in \gamma$ includes an isometric copy of $\mathbb{S} ^{\infty}$.

\end{thm*}

This result is traditionally stated in terms of \emph{oscillation stability}, a concept  coming
from Banach space theory. However, it turns out that it can also be stated thanks to a new concept
of oscillation stability due to Kechris, Pestov and Todorcevic introduced in \cite{KPT} and more
fully developed in \cite{Pe1}. The global formulation of this notion takes place at a very general
level and allows to capture various phenomena coming from combinatorics and functional analysis.
Nevertheless, it can be presented quite simply in the realm of complete separable ultrahomogeneous
metric spaces, where it coincides with the Ramsey-theoretic concept of approximate indivisibility.
Recall that a metric space $\m{X}$ is \emph{ultrahomogeneous} when every isometry between finite
metric subspaces of $\m{X}$ can be extended to an isometry of $\m{X}$ onto itself. Now, for
$\varepsilon \geqslant 0$, call a metric space $\m{X}$ \emph{$\varepsilon$-indivisible} when for
every strictly positive $k \in \omega$ and every $\chi : \funct{\m{X}}{k}$, there is $i<k$ and
$\mc{X} \subset \m{X}$ isometric to $\m{X}$ such that
\[\mc{X} \subset (\overleftarrow{\chi} \{ i \})_{\varepsilon}.\]
Then $\m{X}$ is \emph{approximately indivisible} when $\m{X}$ is $\varepsilon$-indivisible for
every $\varepsilon > 0$, and $\m{X}$ is \emph{indivisible} when $\m{X}$ is $0$-indivisible. For
example, in this terminology, the aforementioned  theorem of Odell and Schlumprecht asserts that
$\mathbb{S}^{\infty}$ is not approximately indivisible. However, in spite of this solution, it is
sometimes felt that something essential is still to be discovered about the metric structure of
$\mathbb{S} ^{\infty}$. Indeed, quite surprisingly, the proof leading to the solution is not based
on an analysis of the intrinsic geometry of $\ell _2$. This fact is one of the motivations for the
present note: In this article, hoping that a better understanding of $\mathbb{S} ^{\infty}$ might
be hidden behind a general approach of approximate indivisibility, we study the approximate
indivisibility problem for another complete, separable ultrahomogeneous metric space, the
\emph{Urysohn sphere} $\s$,   defined as follows: Up to isometry, it is the unique metric space to
which every sphere of radius $1/2$ in the Urysohn space $\Ur$ is isometric. Equivalently, it is, up
to isometry, the unique  complete separable ultrahomogeneous metric space with diameter $1$ into
which every separable metric space with diameter less or equal to $1$ embeds isometrically. In this note, we try
to answer the following question implicitly present in \cite{KPT} and explicitly stated in \cite{Hj} and 
\cite{Pe1}:

\textbf{Question.} Is the Urysohn sphere $\s$ oscillation stable? That is, given a finite
partition $\gamma$ of $\s$ and $\varepsilon > 0$, is there $A \in \gamma$ such that
$(A)_{\varepsilon}$ includes an isometric copy of $\s$?

Our approach here is combinatorial and follows the general intuition according to which the
structure of complete  separable ultrahomogeneous metric spaces can be approached via combinatorial
means. This intuition is based on two ideas. The first one is that the combinatorial point of view
is relevant for the study of countable ultrahomogeneous metric spaces in general. This idea is
already central in the work of Fra\"iss\'{e} completed in the fifties, even though Fra\"iss\'{e} theory
takes place at the level of relational structures and includes much more than metric spaces (for a reference on Fra\"iss\'{e} theory, see \cite{Fr}). More
recently, it was also rediscovered by Bogatyi in a purely metric context, see \cite{Bo0} and
\cite{Bo1}. The second idea is that the complete separable ultrahomogeneous metric spaces are
closely linked to the countable ultrahomogeneous metric spaces. This connection also appears in
Bogatyi's work but is on the other hand supported by the following result (which answers a question posed in \cite{Bo1}):

\begin{thm}

\label{thm:countable dense ultrahomogeneous}

Every complete separable ultrahomogeneous metric space $\m{Y}$ includes a countable ultrahomogeneous dense metric subspace.

\end{thm}

For example, consider the \emph{rational Urysohn space} $\Ur _{\Q}$ which can be defined up to
isometry as the unique countable  ultrahomogeneous metric space with rational distances for which
every   countable metric space with rational distances  embeds isometrically. The Urysohn
space $\Ur$ arises then as the completion of $\Ur _{\Q}$, a fact which is actually essential as it is at
the heart of several important contributions about $\Ur$. In particular, in the original article
\cite{U} of Urysohn, the space $\Ur$ is precisely constructed as the completion of $\Ur _{\Q}$
which is in turn constructed by hand.

Similarly, the Urysohn sphere $\s$ arises as the completion of the so-called \emph{rational Urysohn
sphere} $\s _{\Q}$, defined up  to isometry as the unique countable ultrahomogeneous metric space
with distances in $\Q \cap [0,1]$ into which every at most countable metric space with distances in
$\Q \cap [0,1]$ embeds isometrically.

At first glance, such a representation is relevant with respect to the oscillation stability
problem  for complete separable ultrahomogeneous metric spaces because it provides a direct way to
transfer an approximate indivisibility problem to an exact indivisibility problem. For example, in
the present case, it naturally leads to the question (explicitly stated in \cite{N1} and in
\cite{Pe1}) of knowing whether $\s _{\Q}$ is indivisible, a question which was answered recently by
to Delhomm\'{e}, Laflamme, Pouzet and Sauer in \cite{DLPS}, where a   detailed analysis of metric
indivisibility is provided and several obstructions to indivisibility are isolated. Cardinality is
such an obstruction: A classical result in topology asserts that as soon as a metric space $\m{X}$
is uncountable, there is a partition of $\m{X}$ into two pieces such that none of the pieces
includes a copy of the space via a continuous $1-1$ map. Unboundedness is another example: If a
metric space $\m{X}$ is indivisible, then its distance set is bounded. Now, it turns out that $\s
_{\Q}$ avoids these obstacles but encounters a third one: For a metric space $\m{X}$, $x \in
\m{X}$, and $\varepsilon > 0$, let $\lambda _{\varepsilon} (x)$ be the supremum of all reals $l
\leqslant 1$ such that there is an $\varepsilon$-chain $(x_i)_{i \leqslant n}$ containing $x$ and
such that $d^{\m{X}}(x_0 , x_n) \geqslant l$. Then, define
\[\lambda (x) = \inf \{ \lambda_\vep(x) :  \varepsilon > 0 \}.\]
\begin{thm*}[Delhomm\'{e}-Laflamme-Pouzet-Sauer \cite{DLPS}]
\label{thm:lambda divisible} Let $\m{X}$ be a countable metric space. Assume that there is $x_0 \in
\m{X}$ such that $\lambda (x_0) > 0$. Then $\m{X}$ is not indivisible.
\end{thm*}
Now, for $\s_{\Q}$, it is easy to see that ultrahomogeneity  together with the fact that the
distance set contains $0$ as an accumulation point imply that every point $x$ in $\s _{\Q}$ is such
that $\lambda (x) = 1$. It follows that:

\begin{cor*}[Delhomm\'{e}-Laflamme-Pouzet-Sauer \cite{DLPS}]

\label{thm:s_Q divisible}

$\s _{\Q}$ is divisible.

\end{cor*}

This result put an end to the first attempt to solve the oscillation  stability problem for $\s$.
Indeed, had $\s _{\Q}$ been indivisible, $\s$ would have been oscillation stable. But in the
present case, the coloring which is used to divide $\s _{\Q}$ does not lead to any conclusion and
the oscillation stability problem for $\s$ has to be attacked from another direction.

Here, following with the intuition that approximate indivisibility for $\s$ can be attacked via the
study of the exact  indivisibility of simpler spaces, we can show:

\begin{thm}

\label{thm:s 1/6 os}

$\s$ is $1/3$-indivisible.
\end{thm}

This result is obtained after having shown that the problem of approximate indivisibility for $\s$
can be reduced to a purely combinatorial problem involving a family $(\s_m)_{m \geqslant 1}$ of
countable metric spaces which in some sense approximate the space $\s$. For $m \in \omega$ strictly positive, set
$$[0,1] _m := \{ \frac{k}m: k \in \{ 0, \ldots , m \} \}.$$
Then $\s_m$ is defined as follows: Up to isometry it is the unique countable ultrahomogeneous
metric space with distances in $[0,1] _m$ into which every countable metric space with distances in
$[0,1] _m$ embeds isometrically. Then:

\begin{thm}

\label{thm:TFAE S mos1}

The following are equivalent:

\begin{enumerate}
\item  $\s$ is oscillation stable (equivalently, approximately indivisible).
\item For every strictly positive $m \in \omega$, $\s _m$ is $1/m$-indivisible.
\item For every strictly positive $m \in \omega$, $\s _m$ is indivisible.
\end{enumerate}
\end{thm}

The paper is organized as follows: In section \ref{section:Discretization}, we introduce the spaces
$\s _m$ and study  their relationship with $\s$. In particular, this leads us to a stronger version
of Theorem \ref{thm:TFAE S mos1}. In section \ref{section:RandB}, we follow the different
directions suggested by Theorem \ref{thm:TFAE S mos1} and study the indivisibility as well as the
$1/m$-indivisibility properties of the spaces $\s _m$. We then show how these results can be used
to derive Theorem \ref{thm:s 1/6 os}. Finally, we close with a short section including some remarks
about possible further studies while an Appendix provides a proof of Theorem \ref{thm:countable
dense ultrahomogeneous}.

\

\emph{Note}: Item (iii) of Theorem \ref{thm:TFAE S mos1} was recently proved by the N. W. Sauer and the second author. The Urysohn sphere is therefore oscillation stable. 

\section{Discretization.}

\label{section:Discretization}

The purpose of this section is to prove Theorem \ref{thm:TFAE S mos1} and therefore  to show that
despite the unsuccessful attempt realized with $\s _{\Q}$, the oscillation stability problem for
$\s$ can indeed be understood via the study of the exact indivisibility of simpler spaces. The
starting point of our construction consists in the observation that $\s _{\Q}$ is  the first
natural candidate because it is a very good countable approximation of $\s$, but this good
approximation is paradoxically responsible for the divisibility of $\s _{\Q}$. In particular, the
distance set of $\s _{\Q}$ is too rich and allows to create a dividing coloring. A natural attempt
at that point is consequently to replace $\s _{\Q}$ by another space with a simpler distance set
but still allowing to approximate $\s$ in a reasonable sense. In this perspective, general
Fra\"iss\'{e} theory provides a whole family of candidates. Indeed, recall that for a strictly positive $m \in
\omega$, $[0,1] _m$ denotes the set $\{ k/m : k \in \{ 0, \ldots , m \} \}$. Then one can prove that
there is a countable ultrahomogeneous metric space $\s _m$ with distances in $[0,1]_m$ into which
every countable metric space with distances in $[0,1] _m$ embeds isometrically and that those
properties actually characterize $\s _m$ up to isometry. In other words, the spaces $\s _m$ are
really the analogs of $\s _{\Q}$ after having discretized the distance set $\Q \cap [0,1]$ with
$[0,1]_m$. The intuition is then that in some sense, this should allow them to approximate $\s$.
This intuition turns out to be right, as shown by the following proposition whose proof is
postponed to subsection \ref{subsection:Proof of proposition s_m 1/m dense in s}:

\begin{prop}

\label{thm:s_m 1/m dense in s}

For every strictly positive $m \in \omega$, there is an isometric copy $\widetilde{\s _m} $ of $\s _m$ inside $\s$ such that $(\widetilde{\s _m} )_{1/m} = \s$.

\end{prop}

The spaces $\s _m$ consequently appear as good candidates towards a discretization of the
oscillation stability problem for $\s$. However,  it turns out that Proposition \ref{thm:s_m 1/m
dense in s} is not of any help towards a proof of Theorem \ref{thm:TFAE S mos1}. For example,
Proposition \ref{thm:s_m 1/m dense in s} does not imply alone that if for some strictly positive $m
\in \omega$, $\s _m$ is indivisible, then $\s$ is $1/m$-indivisible: Assume that $\chi :
\funct{\s}{k}$. $\chi$ induces a coloring of $\s_m$ so by indivisibility of $\s_m$ there is $\mc{S}
_m \subset \s _m$ isometric to $\s _m$ on which $\chi$ is constant. But how does that allow to
obtain a copy of $\s$? For example, are we sure that $(\mc{S}_m)_{1/m}$ includes a copy of $\s$? We
are not able to answer this question, but recent results of J. Melleray in \cite{Me3} strongly suggest
that $(\mc{S} _m)_{1/m}$ really depends on the copy $\mc{S} _m$ and can be extremely small. In
particular, it may not include a copy of $\s$. Thus, to our knowledge, Proposition \ref{thm:s_m 1/m
dense in s} does not say anything about the oscillation stability of $\s$, except maybe that the spaces
$\s _m$'s are not totally irrelevant for our purposes.

Fortunately, the spaces $\s _m$ do allow to go much further than Proposition \ref{thm:s_m 1/m dense
in s} and are indeed relevant objects. In particular, they allow to reach the following
equivalence, extending   Theorem \ref{thm:TFAE S mos1}:

\begin{thm}

\label{thm:TFAE S mos}

The following are equivalent:
\begin{enumerate}
\item   $\s$ is oscillation stable.
\item $\s _{\Q}$ is approximately indivisible.
\item For every strictly positive $m \in \omega$, $\s _m$ is $1/m$-indivisible.
\item For every strictly positive $m \in \omega$, $\s _m$ is indivisible.
\end{enumerate}
\end{thm}

Subsections \ref{subsubsection:From  oscillation stability of S to approximate indivisibility of
S_Q.} to  \ref{subsubsection:From indivisibility of S_m to  oscillation stability of S.} are
devoted to the proof of this result. But before going deeper into the technical details, let us
mention here that part of our hope towards the discretization strategy comes from the proof of a
famous result in Banach space theory, namely Gowers' stabilization theorem for $c_0$ \cite{Gow},
where combinatorial Ramsey-type theorems for the spaces $\mathrm{FIN}^{\pm}_k$ and $\mathrm{FIN}_k$
imply that the unit sphere $\mathbb{S}_{c_0}$ of $c_0$ and its positive part $\mathbb{S}_{c_0} ^+$
are approximately indivisible.

\subsection{Proof of proposition \ref{thm:s_m 1/m dense in s}.}

\label{subsection:Proof of proposition s_m 1/m dense in s}

We start with a definition:  Given a metric space $\m{X} = (X, d^{\m{X}})$, a map $f :
\funct{X}{]0,+\infty[}$ is \emph{Kat\u{e}tov over $\m{X}$}  when: \[\forall x, y \in X, \ \ |f(x) - f(y)| \leqslant d^{\m{X}} (x,y) \leqslant f(x) + f(y).\] Equivalently, one can extend the metric $d^{\m{X}}$ on $X \overset{.}{\cup} \{ f \}$ by defining, for every $x, y$ in X, $ \widehat{d^{\m{X}}} (x, f) = f(x)$ and $\widehat{d^{\m{X}}} (x, y) =
d^{\m{X}} (x, y)$. The corresponding metric space is then written $\m{X} \overset{.}{\cup} \{ f
\}$. Here, the concept of Kat\u{e}tov map is relevant because of the following standard reformulation of the
notion of ultrahomogeneity:

\begin{lemma}
\label{prop:extension} Let $\m{X}$ be a countable metric space. Then $\m{X}$ is ultrahomogeneous
iff for every finite subspace  $\m{F} \subset \m{X}$ and every Kat\u{e}tov map $f$ over $\m{F}$, if
$\m{F} \overset{.}{\cup} \{ f \}$ embeds into $\m{X}$, then there is $y \in \m{X}$ such that for
every $x \in \m{F}$, $d^{\m{X}}(x,y) = f(x)$.
\end{lemma}

This result will be used constantly throughout the proof. Now, some notation: For $m \in \omega$
strictly positive, recall that $[0,1] _m = \{ k/m : k \in \{ 0, \ldots , m \} \}.$ For $\alpha \in
[0,1]$, set also
\[\left\lceil \alpha \right\rceil _m = \min ([\alpha ,1] \cap [0,1] _m)=\frac{\lceil m\alpha \rceil}{ m} ,\]
where $\lceil x \rceil=\min([x, \infty[ \cap \mathbb Z)$ is the ceiling function.
Since $\s$ is the metric completion of $\s _{\Q}$, it is enough to show that for every strictly
positive $m \in \omega$, there is an isometric copy $\widetilde{\s} _m$ of $\s _m$ inside $\s
_{\Q}$ such that $(\widetilde{\s} _m)_{1/m} = \s _{\Q}$. This is achieved thanks to a back and
forth argument.  The following is the main idea.
\begin{claim}
Suppose that $ X\subset \s _{\Q}$ is finite and embeddable in $\s_m$, and let $y\in \s_\Q\smallsetminus X$. Then the mapping $f=f_{X,y,m}:X\cup \{y\}\to ]0,\infty[$ defined by
$f(x)=\lceil d^{\s_\Q}(x,y)\rceil_m$ if $x\in X$ and $f(y)=\max\{ \lceil d^{\s_\Q}(x,y)\rceil_m- d^{\s_\Q}(x,y) \, : \, x\in X  \}$
is  Kat\u{e}tov.
\end{claim}
Assume  this claim is true.  Fix $ (x_n) _{n \in \omega}$ an enumeration of $\s _m$ and $ (y_n) _{n \in \omega}$
an enumeration of $\s _{\Q}$.
We are going to construct  $\sigma:\omega \to \omega$  together with a set 
$\widetilde{\s}_m =\{ \widetilde x _{\sigma (n)} : n \in \omega \} \subset \s_\Q$ so that:
\begin{enumerate}
\item $\sigma$ is a bijection.  
\item $\widetilde x_{\sigma(n)} \mapsto x_{\sigma(n)}$ defines an isometry.

\item For every $n \in \omega$, $\{ y_i : i \leqslant n\}\subset ( \{ \widetilde x_i : i \leqslant 2n+1\})_{1/m}$.
\end{enumerate}
Observe that, since $\sigma$ is a permutation,  (i) and (ii) guarantee that $\widetilde x_n \mapsto x_n $ defines a surjective isometry between $\widetilde{\s}_m $ and $\s_m$. On the other hand, (iii) guarantees that $(\widetilde{\s}_m )_{1/m}=\s_\Q$.

Let $\sigma(0)=0$, $\widetilde{x}_0=y_0$. Suppose now all data up to $2n$ already defined in the appropriate way, i.e. fulfilling the obvious partial versions of (i), (ii) and (iii).
Let
\[\sigma (2n+1) = \min (\omega \smallsetminus \{ \sigma (i) : 0 \leqslant i \leqslant 2n \}).\]
Set also $\tilde{x} _{\sigma (2n+1)} \in \s _{\Q}$ such that:
\[\forall i \in \{ 0, \ldots, 2n \}, \ d ^{\s _{\Q}}(\tilde{x} _{\sigma (i)} , \tilde{x} _{\sigma
(2n+1)}) = d ^{\s _m}(x_{\sigma (i)} , x_{\sigma (2n+1)}).\]
Next, if $y_n\in \left(\{\widetilde{x}_{\sigma(i)}: i\leqslant 2n+1\}\right)_{1/m}$, then we define $\sigma(2n+2)$ and $\widetilde x_{\sigma(2n+2)}$ as we did for $2n+1$.  Otherwise, let  $f$ be the Kat\u{e}tov map given by the previous claim when applied to $X= \{ \tilde{x} _{\sigma (i)} : 0 \leqslant i\leqslant 2n+1 \}$ and $y_{n}$.   Let $\tilde{x} \in \s _{\Q}
$ realizing $f$. Now observe that  the map $g$ defined on $\{  {x} _{\sigma (i)} : 0 \leqslant i\leqslant 2n+1 \}$  by $g(x_{\sigma(i)})=f(\widetilde x_{\sigma(i)})$ is  Kat\u{e}tov with values in $[0,1]_m$, so
\[\sigma (2n+2) = \min \{ k \in \omega : \forall i \in \{ 0, \ldots, 2n+1 \}, \  d ^{\s _m }(x_{\sigma
(i)} , x_k) = g ( x _{\sigma (i)} ) \}\]
is well defined and we set $\tilde{x} _{\sigma (2n+2)}=\tilde{x}$. 

We now turn to the proof of the claim. Fix $x,x'\in X$. We have to prove:
\begin{align}
 \left| f (x) - f (x') \right| &\leqslant  d ^{\s
_{\Q}}(x,x') \leqslant f(x) + f(x')
  \label{jdigfdifj1} \\
\left| f (x) - f(y) \right| & \leqslant   d ^{\s
_{\Q}}(x,y ) \leqslant  f (x) + f(y)  \label{jdigfdifj2}
\end{align}
For \eqref{jdigfdifj1}: The right inequality is not a problem:
\[d ^{\s _{\Q}}(x,x') \leqslant d ^{\s
_{\Q}}(x , y ) + d ^{\s _{\Q}}(y,  x') \leqslant f(x) + f(x').\]  For
the left inequality, we use the following simple fact:
\[\forall \alpha , \beta \in \R, \ \forall p \in \omega, \ \left| \beta - \alpha \right|
 \leqslant \frac{p}{m} \longrightarrow \left| \left\lceil \beta \right\rceil _m - \left\lceil \alpha \right\rceil _m \right| \leqslant \frac{p}{m}.\]
Indeed, assume that $\left| \beta - \alpha \right| \leqslant p/m$. We want $\left| \left\lceil m
\beta \right\rceil - \left\lceil m \alpha \right\rceil \right| \leqslant p$. Without loss of
generality, $ \alpha \leqslant \beta$. Then $0 \leqslant \left\lceil m \beta \right\rceil -
\left\lceil m \alpha \right\rceil < m \beta + 1 - m \alpha \leqslant p + 1$, so $\left| \left\lceil
m \beta \right\rceil - \left\lceil m \alpha \right\rceil \right| \leqslant p$ and we are done. In
our case, that property is useful because then the left inequality directly follows from
\[\left|d ^{\s _{\Q}}(x , y  ) - d ^{\s _{\Q}}(y  ,
x' ) \right| \leqslant d ^{\s _{\Q}}(x , x') \in [0,1] _m,\]
because $X$ is embeddable in $\s_m$.
 For \eqref{jdigfdifj2}:
\[\left|f (x) - f(y) \right| = f(x) - f (y  ).\]
This is because $f(x) \geqslant 1/m$ and $ 0 \leqslant f (y ) < 1/m$. Furthermore, by definition of $f$, \[ f(y) \geqslant f(x) -  d ^{\s _{\Q}}(x,y).\]
So the left inequality is satisfied. For the right inequality, simply observe that
\begin{center}
\mbox{ } \hfill $d ^{\s _{\Q}}(x,y) \leqslant
f(x)$. \hspace{\stretch{1}}\qed
\end{center}

\subsection{From oscillation stability of $\s$ to approximate indivisibility of $\s _{\Q}$.}
\label{subsubsection:From  oscillation stability of S to approximate indivisibility of S_Q.}

The purpose of what follows is to prove the implication $(i) \rightarrow (ii)$  of Theorem
\ref{thm:TFAE S mos} stating that if $\s$ is oscillation stable, then $\s _{\Q}$ is approximately
indivisible. This is done thanks to the following result:

\begin{prop}

\label{thm: s_Q in s}
Suppose that $\s_\Q^0$ and $\s_\Q^1$ are two copies of $\s_\Q$ in $\s$ such that $\s_\Q^0$ is dense in $\s$. Then for every $\varepsilon>0$ the  subspace $\s_\Q^0\cap (\s_\Q^1)_\vep$ includes a copy of $\s_\Q$.
\end{prop}
\begin{proof}
We construct the required copy of $\s _{\Q}$ inductively. Let $\{ y_n : n \in \omega \}$ enumerate
 $\s _{\Q}^1$.  For $k \in \omega$, set \[\delta _k = \frac{\varepsilon}{2} \sum
_{i = 0} ^k \frac{1}{2^i}.\] Set also
 \[\eta _k = \frac{\varepsilon}{3} \frac{1}{2^{k+1}}.\]
  $\s
_{\Q}^0$ being dense in $\s$, choose $z_0 \in \s _{\Q}^0$ such that $d^{\s}(y_0 , z_0) < \delta _0$.
Assume now that $z_0,\dots, z_n \in \s _{\Q}^0$ were constructed such that for every $k, l \leqslant
n$
\begin{displaymath}
\left \{ \begin{array}{l}
 d^{\s}(z_k , z_l)=d^{\s} (y_k , y_l) \\
 d^{\s} (z_k , y_k) < \delta _k.
 \end{array} \right.
\end{displaymath}
Again by denseness of $\s_{\Q}^0$ in $\s$, fix $z \in \s _{\Q}^0$ such that

\begin{center}
$d^{\s} (z,y_{n+1}) < \eta _{n+1}$.
\end{center}
Then for every $k \leqslant n$,
\begin{align*}
\left| d^{\s} (z , z_k) - d^{\s} (y_{n+1} , y_k) \right| & =  \left| d^{\s} (z , z_k) - d^{\s}(z_k
, y_{n+1})  + d^{\s}(z_k , y_{n+1})
 - d^{\s} (y_{n+1} , y_k) \right| \\
 & \leqslant  d^{\s} (z , y_{n+1}) + d^{\s} (z_k , y_k) \\
 & <  \eta _{n+1} + \delta _k \\
 & <  \eta _{n+1} + \delta_n.
\end{align*}
It follows that there is $z_{n+1} \in \s _{\Q}^0$ such that
\begin{displaymath}
\left \{ \begin{array}{l}
 \forall k \leqslant n \ \ d^{\s} (z_{n+1} , z_{k}) = d^{\s} (y_{n+1} , y_k) \\
 d^{\s} (z_{n+1} , z) < \eta _{n+1} + \delta _n.
 \end{array} \right.
\end{displaymath}
Indeed, consider the map $f$ defined on $\{ z_k : k \leqslant n \} \cup \{ z \}$ by:
\begin{displaymath}
\left \{ \begin{array}{l}
 \forall k \leqslant n \ \ f(z_k) = d^{\s} (y_{n+1} , y_k) \\
 f(z) =  \left| d^{\s} (z , z_k) - d^{\s} (y_{n+1} , y_k) \right|.
 \end{array} \right.
\end{displaymath}
Then $f$ is Kat\u{e}tov over the subspace of $\s _{\Q}^0$ supported by $\{ z_k : k \leqslant n \}
\cup \{ z \}$, so simply take $z_{n+1} \in \s _{\Q}^0$  realizing it. Observe then that
\begin{eqnarray*}
d^{\s}(z_{n+1} , y_{n+1}) & \leqslant & d^{\s}(z_{n+1},z) + d^{\s}(z,y_{n+1})\\
& < & \eta _{n+1} + \delta _n + \eta _{n+1}\\
& < & \delta _{n+1}.
\end{eqnarray*}
After $\omega$ steps, we are left with $\{ z_n : n \in \omega \} \subset \s _{\Q}^0 \cap
(\s_\Q^1)_{\varepsilon}$ isometric to $\s _{\Q}$. \end{proof}

We now show how to deduce $(i) \rightarrow (ii)$ of Theorem \ref{thm:TFAE S mos} from Proposition
\ref{thm: s_Q in s}:  Let $\varepsilon > 0$, $k \in \omega$ strictly positive and $\chi : \funct{\s
_{\Q}}{k}$. Then in $\s$, seeing $\s_{\Q}$ as a dense subspace:
\[ \s = \bigcup _{i<k} (\overleftarrow{\chi} \{ i \})_{\varepsilon / 2}.\]
By oscillation stability of $\s$, there is $i<k$ and a copy $\mc{S}$ of $\s$ included in $\s$ such
that
\begin{center}
$\mc{S} \subset ((\overleftarrow{\chi} \{ i \})_{\varepsilon/2})_{\varepsilon/4}$.
\end{center}
Since $\widetilde{\s}$ includes copies of $\s_\Q$, and since $\s_\Q$ is dense in $\s$, it follows by Proposition \ref{thm: s_Q in s} that there is a copy $\mc{S} _{\Q}$ of $\s _{\Q}$ in
$\s _{\Q}\cap(\mc{S})_{\varepsilon /4}$. Then in $\s _{\Q}$
\begin{center}
\mbox{ }\hfill$\mbox{ }\hfill \mc{S} _{\Q} \subset (\overleftarrow{\chi} \{ i \})_{\varepsilon}$.
\hspace{\stretch{1}} \qed
\end{center}

\subsection{From approximate indivisibility of $\s _{\Q}$ to $1/m$-indivisibility of $\s _m$.}

Here, we provide a proof for the implication $(ii) \rightarrow (iii)$ of Theorem \ref{thm:TFAE S
mos} according to which if $\s _{\Q}$ is approximately indivisible, then $\s _m$ is
$1/m$-indivisible for every strictly positive $m \in \omega$. This is obtained as the consequence
of the following proposition:

\begin{prop}

\label{thm:Approx indiv s_Q - s_m}

Let $\varepsilon > 0$ and assume that $\s _{\Q}$ is $\varepsilon$-indivisible. Then $\s _m$ is
$1/m$-indivisible  whenever $m \leqslant 1/ \varepsilon$.

\end{prop}

\begin{proof}

Let $\varepsilon > 0 $, assume that $\s _{\Q}$ is $\varepsilon$-indivisible and fix $m \in \omega$
strictly positive such that $\varepsilon \leqslant 1/m$. Define $\left\lceil d^{\s _{\Q}}
\right\rceil _m $ by
\[\forall x, y \in X \ \ \left\lceil d^{\s _{\Q}} \right\rceil _m (x, y) = \left\lceil d^{\s _{\Q}}
(x,y) \right\rceil _m.\]
\begin{claim}
$\left\lceil d^{\s _{\Q}} \right\rceil _m$ is a metric on $\s _{\Q}$.
\end{claim}
\begin{proof}
Since the function $\lceil \cdot \rceil_m$ is subadditive and increasing, it easily follows that the composition $\lceil d^{\s_\Q}\rceil_m =\lceil \cdot \rceil_m \circ d^{\s_\Q}$ is a metric. 
\end{proof}
Let $\m{X} _m$ be the metric space
\[\m{X}_m=(\s _{\Q} , \left\lceil d^{\s _{\Q}} \right\rceil _m),\] and let $\pi _m$  denote the identity map from $\s _{\Q}$
to $\m{X} _m$. Observe that $\m{X} _m$ and $\s _m$ embed into each other, and that consequently,
$1/m$-indivisibility of $\s _m$ is equivalent to $1/m$-indivisibility of $\m{X} _m$. So let $k \in
\omega$ be strictly positive and $\chi : \funct{\m{X} _m}{k}$. Then $\chi$ induces a coloring $\chi
\circ \pi_m : \funct{\s _{\Q}}{k}$. Since $\s _{\Q}$ is $\varepsilon$-indivisible, there is $i<k$ and
a copy $\mc{S} _{\Q}$ of $\s _{\Q}$ inside $\s _{\Q}$ such that

\[\mc{S} _{\Q} \subset (\overleftarrow{\chi \circ \pi_m } \{ i \})_{\varepsilon}.\]
Now, observe that $\pi_m '' \mc{S} _{\Q}$ is a copy of $\m{X}_m$ inside $\m{X}_m$. Furthermore, note that
\[\forall x \neq y \in \s _{\Q} \text{ if }d^{\s _{\Q}} (x , y) \leqslant \frac1m \text{ then }d^{\m{X} _m}
(\pi_m (x), \pi_m(y)) = \frac1m.\] Since $\varepsilon \leqslant 1/m$, it follows that
\[\pi_m '' (\overleftarrow{\chi \circ \pi_m } \{ i \})_{\varepsilon} \subset (\overleftarrow{\chi} \{ i
\})_{1/m}.\] And so
\[\pi_m '' \mc{S} _{\Q} \subset (\overleftarrow{\chi} \{ i \})_{1/m}.\qedhere\]
 \end{proof}

\subsection{From $1/2(m^2+m)$-indivisibility of $\s _{2(m^2 + m)}$ to indivisibility of $\s _m$.}

We now turn to the proof of the implication  $(iii) \rightarrow (iv)$ of Theorem \ref{thm:TFAE S
mos} stating that if for every strictly positive $m \in \omega$, $\s _m$ is $1/m$-indivisible, then
for every strictly positive $m \in \omega$, $\s _m$ is indivisible. This is done via the following proposition:

\begin{prop}
Suppose that for some strictly positive integer $m$,  $\s _{2(m^2 + m)}$ is ${1}/2(m^2 +
m)$-indivisible. Then $\s _m$ is indivisible.

\end{prop}

\begin{proof}

Let $m \in \omega$ be strictly positive and such that $\s _{2(m^2 + m)}$ is ${1}/{2(m^2 +
m)}$-indivisible. We are going to create a metric space $\m{W}$ with distances in $[0,1]_m$ and a
bijection $\pi : \funct{\s _{2(m^2 + m)}}{\m{W}}$ such that for every subspace $\m{Y}$ of $\s
_{2(m^2 + m)}$, if $(\m{Y})_{1/2(m^2 + m)}$ includes a copy of $\s _m$, then so does $\pi ''
\m{Y}$.

Assuming that such a space $\m{W}$ is constructed, the result is proved as follows:
 Observe first that $\m{W}$ and $\s _m$ embed into each other. Indivisibility of $\m{W}$ is consequently
  equivalent to indivisibility of $\s _m$ and it is enough to show that $\m{W}$ is indivisible.
   Let $k \in \omega$ be strictly positive and $\chi : \funct{\m{W}}{k}$.
    Then $\chi \circ \pi : \funct{\s _{2(m^2 + m)}}{k}$ and by ${1}/{2(m^2 + m)}$-indivisibility of $\s _{2(m^2 + m)}$,
    there is $i<k$ such that $(\overleftarrow{\chi \circ \pi} \{ i \})_{1/2(m^2 + m)}$ includes a copy of $\s _{2(m^2 + m)}$. Since $\s _m$ embeds into $\s _{2(m^2 + m)}$, $(\overleftarrow{\chi \circ \pi} \{ i \})_{1/2(m^2 + m)}$ also includes a copy of $\s _m$. Thus, $\overleftarrow{\chi} \{ i \} = \pi '' \overleftarrow{\chi \circ \pi} \{ i \}$ includes a copy of $\s _m$, and therefore a copy of $\m{W}$.

We now turn to the construction of $\m{W}$. This space is obtained by modifying the metric on $\s _{2(m^2 + m)}$ to a metric $d$,
 so that $\m{W} = (S_{2(m^2 + m)} , d)$ and $\pi$ is simply the identity map from $\s _{2(m^2 + m)}$ to $\m{W}$.
 The metric  $d$ is defined as follows: consider the map $f : \funct{[0,1]_{2(m^2+m)}}{[0,1]_{m}}$ defined by
$ f(x) = \frac{l}{m}$ where $l$ is the least integer such that \[ x\leqslant l \left( \frac{1}{m} + \frac{1}{m^2+m} \right).\] Observe that $f$ is increasing, that $f(0) = 0$, and that
\[\forall \alpha \in [0,1]_m \ \ \forall \varepsilon \in \{ -2, -1, 0, 1, 2\} \ \ f\left( \alpha + \frac{\varepsilon}{2(m^2+m)} \right) = \alpha .\]
Note also that $f$ is subadditive: Let $x, y, \in [0,1]_{2(m^2+m)}$. Assume that $f(x) = l/m$. Then there is $n \in \{ 1,\ldots , 2m+4 \}$ such that
\[x = \frac{l-1}{m} + \frac{l-1}{m^2+m} + \frac{n}{2(m^2+m)}.\]
Similarly, there are $l' \in \{ 0,\ldots , m\}$ and $n' \in \{ 1,\ldots , 2m+4 \}$ such that
\[y = \frac{l'-1}{m} + \frac{l'-1}{m^2+m} + \frac{n'}{2(m^2+m)}.\]
So
\begin{eqnarray*}
x + y & = & \left(l+l'\right)\left( \frac{1}{m} + \frac{1}{m^2 + m} \right) - 2\left( \frac{1}{m} + \frac{1}{m^2 + m} \right) + \frac{n + n'}{2(m^2+m)} \\
& = & \left(l+l'\right)\left( \frac{1}{m} + \frac{1}{m^2 + m} \right) + \frac{n - (2m+4) + n' - (2m+4)}{2(m^2+m)} \\
& \leqslant & \left(l+l'\right)\left( \frac{1}{m} + \frac{1}{m^2 + m} \right).
\end{eqnarray*}
Therefore, \[ f(x+y) \leqslant \frac{l+l'}{m} = \frac{l}{m} + \frac{l'}{m} = f(x) + f(y). \]
 It follows that the map $d := f \circ d^{\s
_{2(m^2+m)}}$ is a metric taking    values in $[0,1]_m$.  Now to show that $d$ is as
required, it suffices to prove that for every subspace $\m{Y}$ of $\s _{2(m^2 + m)}$, if
$(\m{Y})_{1/2(m^2 + m)}$ includes a copy of $\s _m$, then $\pi '' \m{Y}$ includes a copy of $\s
_m$. So let $\m{Y}$ be a subspace of $\s _{2(m^2 + m)}$ such that $(\m{Y})_{1/2(m^2 + m)}$ includes
a copy $\mc{S} _m$ of $\s _m$. Then for every $x \in \mc{S} _m$, there is an element $\varphi (x)
\in \m{Y}$ such that $d^{\s _{2(m^2 + m)}} (x, \varphi(x)) \leqslant {1}/{2(m^2 + m)}$. Thus,
\[ \forall x \neq y \in \mc{S} _m \ \ \left| d^{\s _{2(m^2 + m)}} ( \varphi(x) , \varphi(y)) - d^{\s _{2(m^2 + m)}} (x, y) \right|
\leqslant \frac{1}{m^2+m}.\] Since $d^{\s _{2(m^2 + m)}} (x, y) \in [0,1]_m$,
\[ f \left( d^{\s _{2(m^2 + m)}} ( \varphi(x) , \varphi(y)) \right) =  d^{\s _{2(m^2 + m)}} (x, y).\]
That is
\[ d(\pi(\varphi(x)) , \pi(\varphi(y))) =  d^{\s _{2(m^2 + m)}} (x, y).\]
Thus, $\pi '' \ran (\varphi) \subset \pi '' \m{Y}$ is isometric to $\s _m$. \end{proof}

\subsection{From indivisibility of $\s _m$ to oscillation stability of $\s$.}

\label{subsubsection:From indivisibility of S_m to oscillation stability of S.}

We are now ready to close the loop of implications of Theorem \ref{thm:TFAE S mos}. In what
follows, we show that if  $\s_m$ is indivisible for  every  strictly positive $m \in \omega$,
then $\s$ is oscillation stable. This is achieved thanks to the following result:

\begin{prop}
\label{cor:if s_m indiv then s 1/2m-indiv} Assume that for some strictly positive $m \in \omega$,
$\s _m$ is indivisible. Then $\s$ is $1/m$-indivisible.
\end{prop}

\begin{proof}

This is obtained by showing that for every strictly positive  $m \in \omega$, there is an isometric
copy $\s _m ^*$ of $\s _m$ inside $\s$ such that for every $\mc{S} _m \subset \s _m ^* $ isometric
to $\s _m$, $(\mc{S} _m)_{1/m}$ includes an isometric copy of $\s$. This property indeed suffices
to prove Proposition \ref{cor:if s_m indiv then s 1/2m-indiv}: Let $\chi : \funct{\s}{k}$ for some
strictly positive $k \in \omega$. $\chi$ induces a $k$-coloring of the copy $\s _m ^*$. By indivisibility of $\s _m$, find $i<k$ and $\mc{S} _m \subset \s _m ^*$
such that $\chi$ is constant on $\mc{S} _m$ with value $i$. But then, in $\s$, $(\mc{S} _m)_{1/m}$
includes a copy of $\s$. So $(\overleftarrow{\chi} \{ i \})_{1/m}$ includes a copy of $\s$.

We now turn to the construction of $\s _m ^*$. The core of the proof is contained in Lemma
\ref{lem:hedgehog} which we present now. Fix an enumeration $\{y_n : n \in \omega \}$ of $\s
_{\Q}$. Also, keeping the notation introduced in the proof of Proposition \ref{thm:Approx indiv s_Q
- s_m}, let $\m{X} _m$ be the metric space $(\s _{\Q} , \left\lceil d^{\s _{\Q}} \right\rceil _m)$.
The underlying set of $\m{X} _m$ is really $\{y_n : n \in \omega \}$ but to avoid confusion, we
write it $\{x_n : n \in \omega \}$, being understood that for every $n \in \omega$, $x_n = y_n$. On
the other hand, remember that $\s _m$ and $\m{X} _m$ embed isometrically into each other.

\begin{lemma}
\label{lem:hedgehog} There is a countable metric space $\m{Z}$ with distances in $[0,1]$ and
including $\m{X} _m$ such that for every strictly increasing $\sigma : \funct{\omega}{\omega}$
 such that  $x_n \mapsto x_{\sigma (n)}$ is an isometry,   $ (\{ x_{\sigma (n)} : n \in \omega \})_{1/m}$
 includes an isometric copy of $\s _{\Q}$.
\end{lemma}
Assuming Lemma \ref{lem:hedgehog}, we now show how we can construct $\s ^* _m$. $\m{Z}$ is
countable with distances in $[0,1]$ so we may assume that it is a subspace of $\s$. Now, take $\s
_m ^*$ a subspace of $\m{X}_m$ and isometric to $\s _m$. We claim that $\s _m ^* $ works: Let
$\mc{S} _m \subset \s _m ^* $ be isometric to $\s _m$. We first show that $(\mc{S}
_m)_{1/m}$ includes a copy of $\s _{\Q}$. The enumeration $\{x_n : n \in \omega \}$
induces a linear ordering $<$ of $\mc{S} _m$ in type $\omega$. According to lemma \ref{lem:hedgehog}, it suffices to
show that $(\mc{S} _m,<)$ includes a copy of $\{x_n : n \in \omega \}_<$ seen as an ordered metric
space. To do that, observe that since $\m{X} _m$ embeds isometrically into $\s _m$, there is a
linear ordering $<^*$ of $\s _m$ in type $\omega$ such that $\{x_n : n \in \omega \}_<$ embeds into
$(\s _m , <^*)$ as ordered metric space. Therefore, it is enough to show:
\begin{claim}
$(\mc{S} _m,<)$ includes a copy of $(\s _m , <^*)$.
\end{claim}
\begin{proof}
Write
\begin{align*}
(\s _m , <^*)& = \{s_n : n \in \omega \}_{<^*}\\
(\mc{S} _m,<) &  = \{t_n : n \in \omega \}_{<}.
\end{align*}
Let $\sigma (0) = 0$. If $\sigma (0) < \dots < \sigma (n)$ are chosen such that $s_k \mapsto
t_{\sigma (k)}$ is a finite isometry, observe that the following set is infinite
\[\{ i \in \omega : \forall k \leqslant n \ \ d^{\s _m}(t_{\sigma (k)} , t_i) = d^{\s _m}(s_k ,
s_{n+1})\}.\]
Therefore, simply take $\sigma (n+1) $ in that set and larger than $\sigma (n)$.
\end{proof} Observe that since the metric completion of $\s _{\Q}$ is $\s$,
the closure of $(\mc{S} _m)_{1/m}$ in $\s$ includes a  copy of $\s$. Hence we are done since
$(\mc{S} _m)_{1/m}$ is closed in $\s$. \end{proof}

We now turn to the proof of lemma \ref{lem:hedgehog}. The strategy is first to provide the set $Z$
where the required metric space $\m{Z}$ is supposed to be based on, and then to argue that the
distance $d^{\m{Z}}$ can be obtained (lemmas \ref{lem:3} to \ref{lem:6}). To construct $Z$, proceed
as follows: For $t \subset \omega$, write $t$ as the strictly increasing enumeration of its
elements:

\begin{center}
$t = \{t_i : i \in |t| \}_<$.
\end{center}
Now, let $T$ be the set of all finite nonempty subsets $t$ of $\omega$ such that $x_n \mapsto
x_{t_n}$ is an isometry between $\{ x_n : n \in |t| \}$ and $\{x_{t_n} : n \in |t|\}$. This
set $T$ is a tree when ordered by end-extension. Let

\begin{center}
$Z = X_m \overset{.}{\cup} T$.
\end{center}
For $z \in Z$, define
\begin{displaymath}
\pi(z) = \left \{ \begin{array}{cl}
 z & \textrm{if $z \in X_m$.} \\
 x_{\max z } & \textrm{if $z \in T$.}
 \end{array} \right.
\end{displaymath}
Now, consider an edge-labelled graph structure on $Z$ by defining   $\delta$ with domain
$\mathrm{dom} (\delta) \subset Z \times Z $  and range $\mathrm{ran} (\delta) \subset
[0,1]$ as follows:
\begin{itemize}
\item If $s, t \in T$, then $(s,t) \in \mathrm{dom}(\delta)$ iff $s$ and $t$ are $<_T$ comparable. In this case,
\[\delta (s,t) = d^{\s _{\Q}} (y_{|s|-1 }, y_{|t|-1}).\]
\item If $x, y \in X_m$, then $(x,y)$ is always in $\mathrm{dom}(\delta)$ and
\[\delta (x,y) = d^{\m{X}_m} (x, y).\]
\item If $t \in T$ and $x \in X_m$, then $(x,s)$ and $(s,x)$ are in $\mathrm{dom}(\delta)$ iff $x = \pi (t)$. In this case
\[\delta (x,s) = \delta (s,x) = \frac1{m}.\]
\end{itemize}
For a branch $b$ of $T$ and $i \in \omega$, let $b(i)$ be the unique element of $b$ with height $i$
in $T$. Observe that $b(i)$ is a ${i+1}$-element subset of $\omega$. Observe also that for every $i,j \in \omega$,
$b(i)$ is connected to $\pi (b(i))$ and $b(j)$, and
\begin{enumerate}
\item $\delta (b(i), \pi (b(i)) = 1/{m}$,
\item $\delta (b(i), b(j)) = d^{\s _{\Q}} (y_i , y_j)$,
\item $\delta(\pi(b(i)),\pi(b(j)))$ is equal to any of the following quantities: 

$d^{\m{X}_m}(x_{\max b(i)},x_{\max b(j)})= d^{\m{X}_m}(x_i,x_j)=  \lceil d^{\s_Q}(y_i,y_j) \rceil_m$.
\end{enumerate}
In particular, if $b$ is a branch of $T$, then $\delta$ induces a metric on $b$ and the map from
$\s _{\Q}$ to $b$ mapping  $y_i$ to $b(i)$ is a surjective isometry. We claim that if we can show
that $\delta$ can be extended to a metric $d^{\m{Z}}$ on $Z$ with distances in $[0,1]$, then lemma
\ref{lem:hedgehog} will be proved. Indeed, let
\[\mc{X} _m = \{ x_{\sigma (n)} : n
\in \omega \} \subset \m{X} _m,\] with $\sigma : \funct{\omega}{\omega}$ strictly increasing and
$x_n \mapsto x_{\sigma (n)}$ distance preserving. See $\mathrm{ran}(\sigma)$ as a branch $b$ of $T$. Then
$(b, d^{\m{Z}}) = (b, \delta)$ is isometric to $\s _{\Q}$ and
\[b \subset (\pi '' b )_{ 1/{m}} = (\mc{X} _m)_{1/{m}} .\]
Our goal now is consequently to show that $\delta$ can be extended to a metric on $Z$ with values
in $[0,1]$. Recall that for $x, y \in Z$, and $n \in \omega$ strictly positive, a path from $x$ to
$y$ of size $n$ as is a finite sequence $\gamma = (z_i)_{i<n}$ such that $z_0 = x$, $z_{n-1} = y$
and for every $i<n-1$,
\[(z_i, z_{i+1}) \in \dom(\delta).\]
For $x, y$ in $Z$, $P(x,y)$ is the set of all paths from $x$ to $y$. If $\gamma = (z_i)_{i<n}$ is
in $P(x,y)$, $ \| \gamma \|$ is defined as:
\[ \| \gamma \| = \sum _{i=0} ^{n-1} \delta (z_i , z_{i+1} ).\]
On the other hand, $\| \gamma \| _{\leqslant 1} $ is defined as:
 \[ \| \gamma \| _{\leqslant 1} =
\min ( \| \gamma \| , 1).\] We are going to see that the required metric can be obtained with
$d^{\m{Z}}$ defined by
\[ d^{\m{Z}}(x,y) = \inf \{ \| \gamma \| _{\leqslant 1} : \gamma \in P(x,y)\} .\]
Equivalently, we are going to show that for every  $(x,y) \in \dom (\delta)$, every path $\gamma$
from $x$ to $y$ is metric, that is:
\begin{equation}
\label{hojthurhgr}\delta (x,y) \leqslant \| \gamma \| _{\leqslant 1}
\end{equation}
Let $x, y \in Z$. Call a path $\gamma$ from $x$ to $y$ \emph{trivial} when $\gamma = (x,y)$ and
\emph{irreducible} when no proper subsequence of $\gamma$ is a non-trivial path from $x$ to $y$.
Finally, say that $\gamma$ is a \emph{cycle} when $(x,y) \in \dom (\delta)$. It should be clear
that to prove that $d^{\m{Z}}$ works, it is enough to show that the previous inequality
\eqref{hojthurhgr} is true for every irreducible cycle. Note that even though $\delta$ takes only
rational values, it might not be the case for $d^{\m{Z}}$. We now turn to the study of the
irreducible cycles in $Z$.

\begin{lemma}

\label{lem:3}

Let $x, y \in T$. Assume that $x$ and $y$ are not $<_T$-comparable. Let $\gamma$ be an irreducible
path from $x$ to $y$ in $T$. Then there is $z \in T$ such that $z <_T x$, $z <_T y$ and $\gamma =
(x,z,y)$.

\end{lemma}

\begin{proof}

Write $\gamma = (z_i)_{i<n+1}$. $z_1$ is connected to $x$ so $z_1$ is $<_T$-comparable with $x$.
We claim that $z_1 <_T x$ : Otherwise, $x <_T z_1$ and every element of $T$ which is
$<_T$-comparable with $z_1$ is also $<_T$-comparable with $x$. In particular, $z_2$ is
$<_T$-comparable with $x$, a contradiction since $z_2$ and $x$ are not connected. We now claim that
$z_1 <_T y$. Indeed, observe that $z_1 <_T z_2$ : Otherwise, $z_2 <_T z_1 <_T x$ so $z_2 <_T x$
contradicting irreducibility. Now, every element of $T$ which is $<_T$-comparable with $z_2$ is
also $<_T$-comparable with $z_1$, so no further element can be added to the path. Hence $z_2 = y$
and we can take $z_1 = z$. \end{proof}

\begin{lemma}

\label{lem:4}

Every non-trivial irreducible cycle in $X_m$ has size $3$.

\end{lemma}

\begin{proof}

Obvious since $\delta $ induces the metric $d^{\m{X}_m}$ on $X_m$. \end{proof}

\begin{lemma}

\label{lem:4'}

Every non-trivial irreducible cycle in $T$ has size $3$ and is included in a branch.

\end{lemma}

\begin{proof}

Let $c = (z_i)_{i<n}$ be a non-trivial irreducible cycle in $T$. We may assume that $z_0 <_T
z_{n-1}$. Now,  observe that every element of $T$ comparable with $z_0$ is also comparable with
$z_{n-1}$. In particular, $z_1$ is such an element. It follows that $n = 3$ and that $z_0, z_1,
z_2$ are in a same branch. \end{proof}

\begin{lemma}

\label{lem:5}

Every irreducible cycle in $Z$ intersecting both $X_m$ and $T$ is supported by a set whose form is one of the following ones.
\begin{center}
\begin{figure}[h]
\includegraphics[scale=0.73]{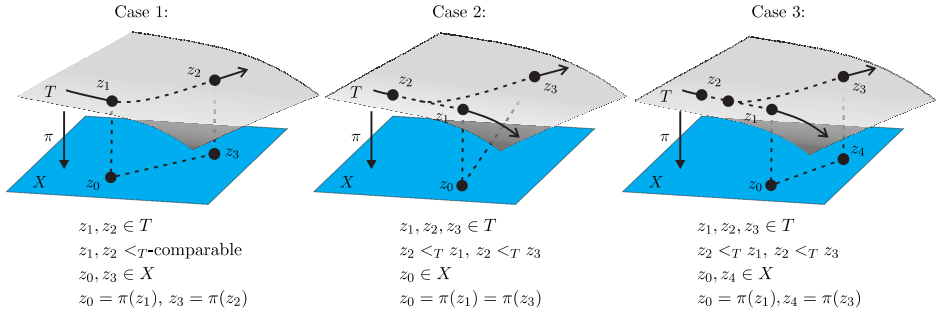}
\caption{Irreducible cycles}\label{figu1}
\end{figure}
\end{center}
%

\end{lemma}

\begin{proof}

Let $C$ be a set supporting an irreducible cycle $c$ intersecting both $X_m$ and $T$. It should be
clear that $|C \cap X_m| \leqslant 2$: Otherwise since any two points in $X_m$ are connected, $c$ would admit a strict subcycle, contradicting irreducibility.

If $C \cap X_m$ has size $1$, let $z_0$ be its unique element. In $c$, $z_0$ is connected to two
elements which we denote $z_1$ and $z_3$. Note that $z_1, z_3 \in T$ so $\pi (z_1) = \pi (z_3) =
z_0$. Since elements in $T$ which are connected never project on a same point, it follows that
$z_1, z_3$ are $<_T$-incomparable. Now, $c$ induces an irreducible path from $z_1$ to $z_3$ in $T$
so from lemma \ref{lem:3}, there is $z_2 \in C$ such that $z_2 <_T z_1$, $z_2 <_T z_3$, and we are
in case 2.

Assume now that $C \cap X_m = \{ z_0 , z_4 \}$. Then there are $z_1, z_3 \in C \cap T$ such that
$\pi(z_1) = z_0$  and $\pi(z_3) = z_4$. Note that since $z_0 \neq z_4$, we must have $z_1 \neq
z_3$. Now, $C \cap T$ induces an irreducible path from $z_1$ to $z_3$ in $T$. By lemma \ref{lem:3},
either $z_1$ and $z_3$ are compatible and in this case, we are in case 1, or $z_1$ and $z_3$ are
$<_T$-incomparable and there is $z_2$ in $C\cap T$ such that $z_2 <_T z_1$, $z_2 <_T z_3$ and we
are in case 3. \end{proof}

\begin{lemma}

\label{lem:6}

Every non-trivial irreducible cycle in $Z$ is metric.

\end{lemma}

\begin{proof}
Let $c$ be an irreducible cycle in $Z$. If $c$ is supported by $X_m$, then by lemma \ref{lem:4} $c$
has size $3$ and  is metric since $\delta$ induces a metric on $X_m$. If $c$ is supported by $T$,
then by lemma \ref{lem:4'} $c$ also has size $3$ and is included in a branch $b$ of $T$. Since
$\delta$ induces a metric on $b$, $c$ is metric. We consequently assume that $c$ intersects both
$X_m$ and $T$. According to lemma \ref{lem:5}, $c$ is supported by a set $C$ whose form is covered
by one of the cases 1, 2 or 3. So to prove the present lemma, it is enough to show every cycle
obtained from a re-indexing of the cycles described in those cases is metric.

Case 1: The required inequalities are obvious after having observed that \[\delta (z_0 , z_3) = \left\lceil \delta (z_1, z_2) \right\rceil _m  \text{ and }  \delta
(z_0 , z_1) = \delta (z_2 , z_3) = \frac1{m}.\]

Case 2: Notice that $\delta (z_0 , z_1) = \delta (z_0 , z_3) = 1/m$. So the inequalities we need to prove are
\begin{align}
\delta (z_1 , z_2) & \leqslant   \delta (z_2, z_3) + \frac2m, \label{lulu1}\\
\delta (z_2 , z_3) & \leqslant  \delta (z_1, z_2) + \frac2m. \label{lulu2}
\end{align}

By symmetry, it suffices to verify that \eqref{lulu1} holds. Observe that since $\pi (z_1) = \pi (z_3) = z_0$, we must have $\left\lceil \delta (z_1 , z_2) \right\rceil_m=\left\lceil \delta(z_2 , z_3)\right\rceil_m$. So: 

\[\delta (z_1 , z_2) \leqslant  \left\lceil \delta (z_1 , z_2) \right\rceil _m = \left\lceil \delta (z_2 , z_3) \right\rceil _m \leqslant  \delta (z_2 , z_3) + \frac2m.\]

Case 3: Observe that $\delta (z_0 , z_1) = \delta (z_3 , z_4) = 1/m$, so the inequalities we need
to prove are
\begin{align}
\delta (z_1 , z_2) & \leqslant   \delta (z_2, z_3) + \delta (z_0 , z_4) + \frac2m, \label{ohrtjuer1}\\
\delta (z_0 , z_4) & \leqslant  \delta (z_1, z_2) + \delta (z_2 , z_3) + \frac2m. \label{ohrtjuer2}
\end{align}
For \eqref{ohrtjuer1}:
\begin{align*}
\delta (z_1 , z_2) & \leqslant  \left\lceil \delta (z_1 , z_2) \right\rceil _m \\
& =  \delta (\pi (z_1) , \pi (z_2) ) \\
& =  \delta (z_0 , \pi (z_2) ) \\
& \leqslant  \delta (z_0 , z_4) + \delta (z_4 , \pi (z_2)) \\
& =  \delta (z_0 , z_4) + \left\lceil \delta (z_3 , z_2) \right\rceil _m \\
& \leqslant  \delta (z_0 , z_4) + \delta (z_2 , z_3) + \frac2m.
\end{align*}
For \eqref{ohrtjuer2}: Write $z_1 = b(j)$, $z_3 = b'(k)$, $z_2 = b(i) = b'(i)$. Then $z_0 = \pi
(z_1) = x_{\max b(j)}$ and $z_4 = \pi (z_3) = x_{\max b'(k)}$. Observe also that $\delta (z_1 ,
z_2) = d^{\s _{\Q}}(y_j , y_i)$ and that $\delta (z_2 , z_3) = d^{\s _{\Q}}(y_i , y_k)$. So:
\begin{align*}
\delta (z_0 , z_4) & =  d^{\m{X}_m}(x_{\max b(j)}, x_{\max b'(k)})\\
& \leqslant  d^{\m{X}_m}(x_{\max b(j)}, x_{\max b(i)}) + d^{\m{X}_m}(x_{\max b'(i)}, x_{\max b'(k)}) \\
& = d^{\m{X}_m}(x_j, x_i) + d^{\m{X}_m}(x_i, x_k) \\
& =  \left\lceil d^{\s _{\Q}}(y_j, y_i)\right\rceil _m + \left\lceil d^{\s _{\Q}}(y_i, y_k) \right\rceil _m \\
& = \left\lceil \delta (z_1, z_2) \right\rceil _m + \left\lceil \delta (z_2 , z_3) \right\rceil _m \\
& \leqslant  \delta (z_1, z_2) + \frac1m + \delta (z_2 , z_3) + \frac1m \\
& = \delta (z_1, z_2) + \delta (z_2 , z_3) + \frac2m. \qedhere
\end{align*}
\end{proof}

\section{Results and bounds.}

\label{section:RandB}

Ideally, the title of this section would have been ``The Urysohn sphere is oscillation stable'' and
we would have ended this article with the proof of one of the different formulations of oscillation
stability for $\s$ presented in Theorem \ref{thm:TFAE S mos}. Unfortunately, so far, our numerous
attempts to reach this goal did not succeed\footnote{The goal has now been achieved by N. W. Sauer and the second author}. This is why this part is entitled ``bounds''. Instead,
what we will be presenting now will show how far we were able to push in the different directions
suggested by Theorem \ref{thm:TFAE S mos}. We start with a summary about the indivisibility
properties of the spaces $\s _m$. 

\subsection{Are the $\s _m$'s indivisible?}

\label{subsection: Are the U_m 's indivisible?} Of course, when $m = 1$, the space $\s _m$ is
indivisible in virtue of the most  elementary pigeonhole principle on $\omega$. The first
non-trivial case is consequently for $m=2$. However, this case is also easy to solve after having
noticed that $\s _2$ is really the Rado graph $\mathcal{R}$ where the distance is $1/2$ between
connected points and $1$ between non-connected distinct points. Therefore, indivisibility for $\s
_2$ is equivalent to indivisibility of $\mathcal{R}$, a problem whose solution is well-known:

\begin{prop}

\label{prop:R indivisible}

The Rado graph $\mathcal{R}$ is indivisible.
\end{prop}

The following case to consider is $\s _3$, which turns out to be another particular case thanks  to
an observation made in \cite{DLPS}. Indeed, $\s _3$ can be encoded by the countable
ultrahomogeneous edge-labelled graph with edges in $\{ 1/3, 1\}$ and forbidding the complete
triangle with labels $1/3, 1/3, 1$. The distance between two points connected by an edge is the
label of the edge while the distance between two points which are not connected is $2/3$. This fact
allows to show:

\begin{thm*}[Delhomm\'{e}-Laflamme-Pouzet-Sauer \cite{DLPS}]

\label{thm:U_3 indivisible}

$\s _3$ is indivisible.

\end{thm*}

Indeed, the proof of this theorem can be deduced from the proof of the indivisibility of the
$\m{K}_n$-free  ultrahomogeneous graph by El-Zahar and Sauer in \cite{EZS1}. We do not write more
here but the interested reader is referred  to \cite{DLPS}, section on the indivisibility of
Urysohn spaces, for more details.

The very first substantial case consequently shows up for $m=4$. Unfortunately, it appears to be so
substantial that so far, we still do not know whether this space is indivisible or not.
Nevertheless, we are able to establish that if this space is not indivisible, then $\s _4$ is quite
exceptional, in a sense that we precise now. We already mentioned that \cite{DLPS} contains an
analysis of indivisibility in the realm of countable metric spaces. It turns out that this study
also led its authors to examine the conditions under which a set of strictly positive reals can be
interpreted as the distance set of a countable universal and ultrahomogeneous metric space:

\begin{defn*}[$4$-values condition]

\label{defn:4values}

Let $S \subset ]0, +\infty[$. $S$ satisfies the $4$-\emph{values condition} when for every $s_0 ,
s_1 , s_0 ', s_1 ' \in S$, if  there is $t \in S$ such that:
\[|s_0 - s_1| \leqslant t \leqslant s_0 + s_1 \text{ and } |s_0 '- s_1 '| \leqslant t \leqslant s_0 '+ s_1 ',\]
then there is $u \in S$ such that:
\[|s_0 - s_0 '| \leqslant u \leqslant s_0 + s_0 ' \text{ and } |s_1 - s_1 '| \leqslant u \leqslant s_1 +
s_1 '.
\]
\end{defn*}

\begin{thm*}[Delhomm\'{e}-Laflamme-Pouzet-Sauer \cite{DLPS}]

\label{thm:4values}

Let $S \subset ]0, +\infty[$. TFAE:
\begin{enumerate}
\item There is a countable ultrahomogeneous metric space $\Ur _S$ with distances in $S$ into which
every countable metric space with distances in $S$ embeds isometrically.
\item $S$ satisfies the $4$-values condition.
\end{enumerate}
\end{thm*}
As detailed in \cite{DLPS}, the $4$-values condition covers a wide variety of examples. For our
purposes, the $4$-values condition is relevant because it allows to establish a list of spaces such
that any space $\Ur _S$ with $S$ finite is in some sense isomorphic to some space in the list. In
particular, it allows to set up a finite list of spaces exhausting all the spaces $\Ur _S$ with $S
\leqslant 4$. More precisely, for finite subsets $S = \{s_0,\ldots , s_m \} _<$, $T = \{t_0,\ldots
, t_n \} _<$ of $]0 , + \infty[$, define $S \sim T$ when $m=n$ and:
\[\forall i, j, k < m, \ \ s_i \leqslant s_j + s_k \leftrightarrow t_i \leqslant t_j + t_k.\]
Observe that when $S \sim T$, $S$ satisfies the $4$-value condition iff $T$ does and in this case,
$S$ and $T$ essentially provide the same metric spaces as it is possible to have $\Ur _S$ and $\Ur
_T$ supported by $\omega$ with the metrics $d^{\Ur _S}$ and $d^{\Ur _T}$ being defined such that:

\[\forall x, y \in \omega, \ \ d^{\Ur _S}(x,y) = s_i \leftrightarrow d^{\Ur _T}(x,y) = t_i.\]
Now, clearly, for a given cardinality $m$ there are only finitely many $\sim$-classes, so we can
find a finite collection $\mathcal{S} _m$ of finite subsets of $]0, \infty[$ of size $m$ such that
for every $T$ of size $m$ satisfying the $4$-value condition, there is $S \in \mathcal{S} _m$ such
that $T \sim S$. For $m \leqslant 3$, examples of such lists can be easily provided. For instance,
one may take:
\begin{align*}
\mathcal{S}_1  & = \{ \{ 1 \} \} \\
\mathcal{S}_2&  = \{ \{ 1, 2 \}, \{ 1, 3\} \} \\
\mathcal{S}_3 & = \{ \{ 2, 3, 4 \}, \{ 1, 2, 3 \}, \{ 1, 2, 5\}, \{ 1, 3, 4\}, \{ 1, 3, 6\}, \{ 1,
3, 7\} \}.
\end{align*}
Notice that in those lists, the set $[0,1]_m$ is represented by the set $\{ 1, 2,\ldots , m\}$. For
$m = 4$, a long and tedious checking procedure of the $4$-values condition allows to find
$\mathcal{S} _m$ but it then contains more than 20 elements so there is no point writing them all
here. Still, it turns out that in most of the cases, we are able to solve the indivisibility
problem for the space $\Ur _S$. Our result can be stated as follows:

\begin{thm}

\label{thm:U_S indiv}

Let $S$ be finite subset of $]0, + \infty [$ of size $|S| \leqslant 4$ and satisfying the
$4$-values condition. Assume that $S \nsim \{ 1, 2, 3, 4\}$. Then $\Ur _S$ is indivisible.

\end{thm}

Due to the number of cases to consider, we do not prove this theorem here but simply mention that
when the proofs are not elementary, three essential ingredients come into play. The first one is
the usual infinite Ramsey theorem, due to Ramsey. The second one is due to El-Zahar and Sauer and
was already mentioned when dealing with $\s _3$. As for the last one, it is due to Milliken (For
more on this theorem and its applications, see \cite{T1}).

The case $S=\{ 1, 2, 3, 4\}$ is consequently the only case with $S = 4$ for which the indivisibility problem
remains unsolved. In the present case, it is a bit ironical as $\{ 1, 2, 3, 4\}$ is precisely the
distance set in which we were interested. So far, the reason for which $\s _4$ stands apart is
still unclear. However, it might be that it is actually the very first case were metricity comes
into play. Indeed, for all the other sets $S$ with $|S| \leqslant 4$, the space $\Ur _S$ can be
coded as an object where the metric aspect does not appear and this is what makes Ramsey, Milliken
and El-Zahar - Sauer theorems helpful. Our feeling is consequently that solving the indivisibility
problem for $\s _4$ requires a new approach. Still, we have to admit that what we are hoping for is
a positive answer and that Theorem \ref{thm:U_S indiv} is undoubtedly responsible for that.

\subsection{$\boldsymbol{1/m}$-indivisibility of the $\s _m$'s.}

We now turn to $1/m$-indivisibility of the spaces $\s _m$. In Theorem \ref{cor:if s_m indiv then s
1/2m-indiv}, we showed how an exact indivisibility result transfers to an approximate one. It turns
out that a slight modification of the proof allows to show:

\begin{prop}

\label{thm:if S_m indiv, then S_3m 1/3m indiv}

Assume that for some strictly positive $m \in \omega$, $\s _m$ is indivisible. Then $\s _{3m}$ is $2/3m$-indivisible.
\end{prop}

\begin{proof}

To prove this theorem, it suffices to show that there is an isometric copy $\s _m ^{**} $ of $\s
_m$ inside $\s _{3m}$ such that for every $\mc{S} _m \subset \s _m ^{**} $ isometric to $\s _m$,
$(\mc{S} _m)_{2/3m}$ includes an isometric copy of $\s _{3m}$. The proof is essentially the same as
the proof of Proposition \ref{cor:if s_m indiv then s 1/2m-indiv} where $\s _m ^{*}$ is
constructed except that instead of the metric space $\m{X} _m = (\s _{\Q} , \left\lceil d^{\s
_{\Q}} \right\rceil_m)$, one works with $(\s _{3m} , \left\lceil d^{\s _{\Q}} \right\rceil_m)$. The
fact that the approximation can be made up to $2/3m$ and not $1/m$ comes from the fact that for
$\alpha \in [0,1]_{3m}$, $\alpha \leqslant \left\lceil \alpha \right\rceil _m \leqslant \alpha +
2/3m$ whereas if $\alpha \in [0,1] \cap \Q$, one only has $\alpha \leqslant \left\lceil \alpha
\right\rceil _m < \alpha + 1/m$. \end{proof}

Thus:

\begin{thm}

\label{thm:s_m, m<10}

For every $m \leqslant 9$, $\s _m$ is $2/m$-indivisible.
\end{thm}


\subsection{Bounds.}

We now turn to the computation of values $\varepsilon$ with respect to which $\s$ is
$\varepsilon$-indivisible. 


\begin{thm*}[Theorem \ref{thm:s 1/6 os}]

$\s$ is $1/3$-indivisible.
\end{thm*}

Note also that if at some point an approximate indivisibility result for $\s _m$ showed up
independently of an exact  one, we would still be able to compute a bound for $\s$:

\begin{prop}
\label{thm:If S _m 1/m-indiv, then S 3/2m indiv} Suppose that for some strictly positive integer
$m$, $\s _m$ is $1/m$-indivisible. Then $\s $ is $\varepsilon$-indivisible for every $\varepsilon
\geqslant 3/2m$.
\end{prop}

\begin{proof}
Let $\varepsilon \geqslant 3/2m$. Consider $\s _m ^*$ constructed in Proposition \ref{cor:if s_m
indiv then s 1/2m-indiv}.  Now, let $k \in \omega$ be strictly positive and $\chi : \funct{\s}{k}$.
$\chi$ induces a coloring of $\s _m ^*$ and $\s _m$ being $1/m$-indivisible, there are $i<k$ and
$\mc{S} _m \subset \s _m ^* $ isometric to $\s _m$ such that $\mc{S} _m \subset
(\overleftarrow{\chi} \{ i \}) _{1/m}$. By construction, $(\mc{S} _m)_{1/2m}$ includes an isometric
copy of $\s$. Now,

\[((\overleftarrow{\chi} \{ i \}) _{1/m}) _{1/2m} \subset (\overleftarrow{\chi} \{ i \}) _{3/2m}
\subset  (\overleftarrow{\chi} \{ i \}) _{\varepsilon}.\] It follows that $(\overleftarrow{\chi} \{
i \}) _{\varepsilon}$ includes an isometric copy of $\s$. \end{proof}


\

\section{Concluding remarks and open problems.}
\label{subsubsection:C_m}

The equivalence provided by Theorem \ref{thm:TFAE S mos} suggests several  lines of future
investigation. Apparently, here is the first and most reasonable question to consider:


\

\textbf{Question.} Is $\s _4$ indivisible? More generally,  is $\s _m$ indivisible for every
strictly positive integer $m$?

\

We finish with two results which might be useful for that purpose. The first one makes a reference to the space $\s _{\Q}$:

\begin{prop}

\label{thm:s_Q for s_m}

Let $m \in \omega$ be strictly positive. Assume that for every strictly positive $k \in \omega$ and
$\chi : \funct{\s _{\Q}}{k}$, there  is a copy $\mc{S} _m$ of $\s _m$ in $\s _{\Q}$ on which $\chi$
is constant. Then $\s _m$ is indivisible.
\end{prop}

\begin{proof}
Once again, we work with $\m{X} _m = (\s _{\Q} , \left\lceil d^{\s _{\Q}} \right\rceil _m)$ and the
identity  map $\pi _m : \funct{\s _{\Q}}{\s _m}$. Think of $\m{X} _m$ as a subspace of $\s _m$.
Now, let $k \in \omega$ be strictly positive and $\chi : \funct{\s _m}{k}$. Then $\chi$ induces a
coloring of $\m{X} _m$, and therefore a coloring $\chi \circ \pi$ of $\s _{\Q}$. By hypothesis,
there is a copy $\mc{S} _m$ of $\s _m$ in $\s _{\Q}$ on which $\chi \circ \pi$ is constant with
value $i<k$. Then $\pi '' \mc{S} _m \subset \overleftarrow{\chi} \{ i \}$. The result follows since
$\pi '' \mc{S} _m$ is isometric to $\s _m$. \end{proof}

The second result provides a space whose indivisibility is equivalent to the indivisibility of $\s
_m$. Let $P$  denote the Cantor space, that is the topological product space $2^{\omega}$. Let
$\mathcal{C}(P)$ denote the set of all continuous maps from $P$ to $\R$ equipped with the $\| .
\|_{\infty}$ norm. Since the work of Banach and Mazur, it is known that $\mathcal{C}(P)$ is a
universal separable metric space. Actually, Sierpinski's proof of that fact allows to show the
following result. For $m \in \omega$ strictly positive, let $\m{C}_m$ denote the space of all
continuous maps from $P$ to $[0,1]_m$ equipped with the distance induced by $\| . \|_{\infty}$.

\begin{prop}
$\m{C}_m$ is a countable metric space and is universal for the class of all countable metric spaces with distances in $[0,1]_m$.
\end{prop}

It follows that $\s _m$ is indivisible iff $\m{C}_m$ is. $\m{C}_m$ being a much more concrete
object than $\s _m$,  studying its indivisibility might be a alternative to solve the
indivisibility problem for $\s _m$.

\section{Appendix - Proof of Theorem \ref{thm:countable dense ultrahomogeneous}.}

\label{section:thm countable dense ultrahomogeneous}

Unlike the rest of this paper, this section does not specifically deal with the oscillation
stability for $\s$ and is simply included here for the sake of completeness. Our purpose is to
prove Theorem \ref{thm:countable dense ultrahomogeneous} by constructing the required subspace of
$\m{Y}$. Let $\m{X} _0 \subset \m{Y}$ be countable and dense. Then, assuming that $\m{X} _n \subset
\m{Y}$ countable has been constructed, get $\m{X} _{n+1}$ as follows: Consider $\mathcal{F}$ the
set of all finite subspaces of $\m{X} _n$. For $\m{F} \in \mathcal{F}$, consider the set $E_n
(\m{F})$ of all Kat\u{e}tov maps $f$ over $\m{F}$ with values in the set $\{ d^{\m{Y}} (x,y) : x, y
\in \m{X} _n \}$ and such that $\m{F} \cup \{ f \}$ embeds into $\m{Y}$. Observe that $\m{X} _n$
being countable, so are $\{ d^{\m{Y}} (x,y) : x, y \in \m{X} _n \}$ and $E_n (\m{F})$. Then, for
$\m{F} \in \mathcal{F}, f \in E_n (\m{F})$, fix $y _{\m{F}} ^f \in \m{Y}$ realizing $f$ over
$\m{F}$. Finally, let $\m{X}_{n+1}$ be the subspace of $\m{Y}$ with underlying set $X_n \cup \{ y_F
^f : \m{F} \in \mathcal{F}, f \in E_n (\m{F}) \}$. After $\omega$ steps, set $\m{X} = \bigcup _{n
\in \omega} \m{X} _n$. $\m{X}$ is clearly a countable dense subspace of $\m{Y}$, and it is
ultrahomogeneous thanks to the equivalent formulation of ultrahomogeneity provided in lemma
\ref{prop:extension}.

A second proof involves logical methods. Fix a countable elementary submodel $M \prec H _{\theta}$
for some  large enough $\theta$ and such that $Y , d^{\m{Y}} \in M$. Let $\m{X} = M \cap \m{Y}$. We
claim that $\m{X}$ has the required property. First, observe that $\m{X}$ is dense inside $\m{Y}$
since by the elementarity of $M$, there is a countable $D \in M$ (and therefore $D \subset M$)
which is a dense subset of $\m{Y}$. For ultrahomogeneity, let $\m{F} \subset \m{X}$ be finite and
let $f$ be a Kat\u{e}tov map over $\m{F}$ such that $\m{F} \cup \{ f \}$ embeds into $\m{X}$.
Observe that $f \in M$. Indeed, $\dom (f) \in M$. On the other hand, let $\mc{F} \cup \{ y \}
\subset \m{X}$ be isometric to $\m{F} \cup \{ f \}$ via an isometry $\varphi$. Then for every $x
\in \m{F}, d^{\m{Y}}(\varphi (x) ,y) \in M$. But $d^{\m{Y}}(\varphi (x) ,y) = f(x)$. Thus, $
\ran(f) \in M$. It follows that $f$ is an element of $M$. Now, by ultrahomogeneity of $\m{Y}$,
there is $y$ in $\m{Y}$ realizing $f$ over $\m{F}$. So by elementarity, there is $x$ in $\m{X}$
realizing $f$ over $\m{F}$.

\end{document}